\documentclass[12pt]{elsart}

\usepackage{amssymb}
\makeatletter

\@addtoreset{equation}{section}

\makeatother

\renewcommand\thefigure{\thesection.\@arabic\c@figure}

\renewcommand\thetable{\thesection.\@arabic\c@table}

\def\reff#1{(\ref{#1})}
\def\sobre#1#2{\lower 1ex \hbox{ $#1 \atop #2 $ } }
\def\supp{{\rm Supp}\,}
\def\basis{{\rm Basis}\,}
\def\life{{\rm Life}\,}
\def\birth{{\rm Birth}\,}
\def\flag{{\rm Flag}\,}
\def\death{{\rm Death}\,}

\def\XX{{\bf X}}

\begin{document}

\def\E{{\mathbb E}}
\def\P{{\mathbb P}}
\def\R{{\mathbb R}}
\def\Z{{\mathbb Z}}
\def\V{{\mathbb V}}
\def\N{{\mathbb N}}
\def\BB{{\mathbb B}}
\def\NN{{\bf N}}
\def\X{{\mathcal X}}
\def\Y{{\bf Y}}
\def\T{{\mathcal T}}
\def\C{{\bf C}}
\def\D{{\bf D}}
\def\G{{\bf G}}
\def\U{{\bf U}}
\def\K{{\bf K}}
\def\H{{\bf H}}
\def\J{{\mathcal J}}
\def\n{{\bf n}}
\def\dd{{\bf d}}
\def\b{{\bf b}}
\def\aa{{\bf a}}
\def\ee{{\rm e}}
\def\g{{\bf g}}
\def\mm{{m}}
\def\ti{{\rm TI}}
\def\sqr{\vcenter{
         \hrule height.1mm
         \hbox{\vrule width.1mm height2.2mm\kern2.18mm\vrule width.1mm}
         \hrule height.1mm}}                  % This is a slimmer sqr.
\def\square{\ifmmode\sqr\else{$\sqr$}\fi}
\def\one{{\bf 1}\hskip-.5mm}
\def\liml{\lim_{L\to\infty}}
\def\given{\ \vert \ }
\def\Given{\ \Big\vert \ }
\def\ze{{\zeta}}
\def\be{{\beta}}
\def\de{{\delta}}
\def\la{{\lambda}}
\def\ga{{\gamma}}
\def\th{{\theta}}
\def\vep{{\varepsilon}}
\def\proof{\noindent{\bf Proof. }}
\def\A{{\bf A}}
\def\B{{\bf B}}
\def\D{{\bf D}}
\def\H{{\bf H}}
\def\h{{\bf h}}
\def\bx{{\bf X}}
\def\bz{{\bf Z}}
\def\bk{{\bf K}}
\def\bF{{\bf F}}
\def\cw{{\mathcal W}}
\def\zero{{\rm 0}}
%\def\MM{{\bf m}}

% New macros
%\def\grains{\underline{h}}
\def\grains{\xi}
\def\VV{{\underline V}}
\def\UU{{\underline U}}
\def\XX{{\underline X}}
\def\FF{{\mathcal F}}
\def\GG{{\mathcal G}}
\def\tends#1{\mathop{\longrightarrow}\limits_{#1}}% Arrow with limits
                                                  % below 
\bibliographystyle{alpha}
\runauthor{Fernández, Ferrari and Garcia}

\begin{frontmatter}
\title{Perfect simulation for interacting point processes,
loss networks and Ising models}
\date{May 18, 2002}
     
\author[USP]{Pablo A. Ferrari}
\author[ROUEN]{Roberto Fern\'andez}
\author[UNICAMP]{Nancy L. Garcia}
\address[USP]{Universidade de S\~{a}o Paulo}
\address[ROUEN]{Universit\'e de Rouen}
\address[UNICAMP]{Universidade Estadual de Campinas}

\begin{abstract} We present a perfect simulation algorithm for measures
  that are absolutely continuous with respect to some Poisson process and can
  be obtained as invariant measures of birth-and-death processes.  Examples
  include area- and perimeter-interacting point processes (with stochastic
  grains), invariant measures of loss networks, and the Ising contour and
  random cluster models.  The algorithm does not involve couplings of the
  process with different initial conditions and it is not tied up to
  monotonicity requirements. Furthermore, it directly provides perfect samples
  of finite windows of the \emph{infinite-volume} measure, subjected to time
  and space ``user-impatience bias''.  The algorithm is based on a two-step
  procedure: (i) a perfect-simulation scheme for a (finite and random)
  relevant portion of a (space-time) marked Poisson processes (free
  birth-and-death process, free loss networks), and (ii) a ``cleaning''
  algorithm that trims out this process according to the interaction rules of
  the target process.  The first step involves the perfect generation of
  ``ancestors'' of a given object, that is of predecessors that may have an
  influence on the birth-rate under the target process.  The second step, and
  hence the whole procedure, is feasible if these ``ancestors'' form a finite
  set with probability one.  We present a sufficiency criteria for this
  condition, based on the absence of infinite clusters for an associated
  (backwards) oriented percolation model.  The criteria is expressed in terms
  of the subcriticality of a majorizing multi-type branching process, whose
  corresponding parameter yields bounds for errors due to space-time
  ``user-impatience bias''.  The approach has previously been used, as an
  alternative to cluster expansion techniques, to extract properties of the
  invariant measures involved.
\end{abstract}


\end{frontmatter}

\section{Introduction}

\emph{Perfect simulations} or \emph{exact sampling} are labels for a recently
developed set of techniques designed to produce output whose distribution is
guaranteed to follow a given probability law.  These techniques are
particularly useful in relation with Markov Chain Monte Carlo, and their range
of applicability is rapidly growing (see {\tt
  http://dimacs.rutgers.edu/\~{}dbwilson/exact}).

There are several techniques for perfect simulation of Markov processes.  The
most popular ones can be classified in two categories: Propp and Wilson's
\emph{Coupling from the Past} (CFTP) and Fill's Interruptible Algorithm.  The
first type of technique applies, in its original version, to invariant
measures of Markov processes with a \emph{finite coalescence time}.  That is,
of processes for which there exists a coupling among trajectories such that
with probability one the trajectories starting from all possible initial
states coalesce in a finite time.  This includes all irreducible Markov
processes with a finite state space.  The coalescence property becomes
difficult to check if the state space is very large.  The problem can be
overcome for processes with the following \emph{monotonicity property}: there
must exist a ``maximal'' and a ``minimal'' state and a coupling such that the
coalescence of coupled trajectories starting from these two states implies the
coalescence of all other trajectories (``monotone coupling'').  Examples of
processes with this property include Glauber dynamics of spin systems with the
FKG property (Propp and Wilson, 1996).  Other perfect simulation techniques
are based on backward coupling of embedded regeneration times (Corcoran and
Tweedie, 2001), finitary coding (H\"{a}ggstr\"{o}m and Steif, 2000; van den
Berg and Steif, 1999), tempering algorithms, cluster representation of a
Markov chain (Cai, 1999 preprint), regenerative construction (Comets,
Fern\'{a}ndez and Ferrari, 2002), embedding the problem into the coloring of a
graph (Fill and Huber, 2000) , representation as an infinite mixture (Hobert
and Robert,2000) among others.

The basic CFTP algorithm, sometimes called \emph{vertical CFTP}, is in
general not applicable to processes with infinite state space.
Indeed, most of them lack uniform ergodicity, a property shown (Foss
and Tweedie, 1998) to be equivalent to the existence of a coalescence
scheme as above.  To cope with this situation, Kendall (1997 and 1998)
introduced \emph{dominated CFTP} (also called \emph{horizontal CFTP}
and \emph{coupling into and from the past}).  This extension also
requires the state space to have a partial order, as well as the
existence of a monotone coupling among the target process and two
reversible \emph{sandwiching processes}, which must be easy to sample.
Algorithms of this type are available for attractive point processes
and, through a minor modification, also for repulsive point processes
(Kendall, 1998). Similarly, H\"aggstr\"om, van Lieshout and M\o ller
(1999) combined ideas from CFTP and the two-component Gibbs sampler to
perfect simulate from process in infinite spaces which do not have
maximal (or minimal) elements.

The interruptible algorithm proposed by Fill (1998) (see also Th\"onnes, 1997)
is an acceptance-rejection scheme which applies to invariant measures of
Markov processes whose \emph{time-reversed} process has a monotonicity
property.  Thus its range of applicability overlaps with that of the CFTP
algorithm at reversible monotone processes like Glauber dynamics of attractive
automata or ferromagnetic spin systems and attractive point processes.  Later
developments have made Fill's algorithm applicable to other processes as well
(Fill, Machida, Murdoch and Rosenthal, 2000).  An important advantage of this
algorithm is that it is free of the so called \emph{impatient-user bias}: no
bias is introduced if the user aborts a long run of the algorithm.

Kendall (1997, 1998) and Kendall and M\o ller (2000) proposed
dominated CFTP schemes applicable to finite-volume measures which are
absolutely continuous with respect to a finite Poisson point process
and that can be obtained as the invariant measure of an interacting
spatial birth-and-death process. These algorithms are based on two
ingredients: (i) the ``thinning'' of a space-time marked Poisson
process, and (ii) the coupled construction of upper and lower
processes whose coalescence signals the output of a perfect sample.
See the recent review of M\o ller (2000) for more references.

In this paper we propose a new perfect-simulation algorithm which
applies to the same type of measures but has the following distinctive
features:

\begin{itemize}
  
\item We sample directly from a time stationary realization of the
  process. There is no coalescence criterion, either between coupled
  realizations or between sandwiching processes.  The scheme neither
  requires nor takes advantage of monotonicity properties.
  
\item The scheme directly samples a finite window of the equilibrium measure
  in \emph{infinite-volume}.  In contrast, Kendall (1997, 1998) focus on
  finite windows with fixed boundary conditions, and the infinite-volume limit
  requires an additional process of ``perfect simulation in space''.  Our
  construction is in the spirit of the algorithms proposed by van den Berg and
  Steif (1999) and by H\"aggstr\"om and Steif (1999) to simulate
  infinite-volume measures for nearest neighbor interactions in a lattice at
  high temperature or ``high noise''.  In a sense, our algorithm is
  complementary to those, because it applies to regimes where they break down
  (e.g.\ at low-temperature).  We point out that before the arrival of the
  perfect simulation wave, Ferrari (1990), van den Berg (1993) and van den
  Berg and Maes (1994) have also proposed construction schemes for
  (infinite-volume) Gibbs measures of spin systems that can be easily
  transcribed into perfect-simulation algorithms.
  
\item The construction has the added value of being a proven theoretical tool
  for the analysis of properties of the target measure.  For instance, in
  Fern\'andez, Ferrari and Garcia (1998 and 2001) we used it to obtain mixing
  properties, finite-volume corrections and the asymptotic (in temperature)
  distribution of ``defects'' of the low-temperature Ising
  translation-invariant extremal measures. Ferrari and Garcia (1998) used a
  similar construction to show ergodicity of a family of loss networks in
  $\R$.

\item More generally, the construction can be used as an alternative to the
  cluster-expansion technology (Brydges 1986, Koteck\`y-Preiss 1986, Dobrushin
  1996) for the study of spin models, at least those with a flipping symmetry.
  In fact, it seems to have a region of validity more extended than usual
  cluster-expansion approaches.

\end{itemize}

Our algorithm does involve the ``thinning'' of a marked Poisson process ---the
\emph{free process}--- which dominates the birth-and-death process, and it
involves a time-backward and a time-forward sweep.  But these procedures are
performed in a form quite different from previous algorithms.  The initial
stage of our construction is done \emph{towards the past}, starting with a
finite window and retrospectively looking to \emph{ancestors}, namely to those
births in the past that could have (had) an influence on the current birth.
The construction of the \emph{clan of ancestors} constitutes the time-backward
sweep of the algorithm.  Once this clan is completely constructed, the
algorithm proceeds in a time-forward fashion ``cleaning up'' successive
generations according to appropriate penalization schemes.  This ``ancestors
approach'' offers some noteworthy advantages:

\begin{itemize}
  
\item[(i)] The algorithm constructs only the portion of the birth-and-death
  process \emph{strictly} needed for the \emph{final} window at $t=0$.  This
  economy has two important consequences: First, we can sample directly from
  the infinite-volume measure, without boundary effects.  Second, the scheme
  works for point processes with quite general grain distribution, for
  instance chosen from an unbounded family of objects.  As a consequence, it
  can be applied to the simulation of loss networks or of Peierls contours of
  the low-temperature Ising model.
  
\item[(ii)] Perfect sampling is assured once the algorithm determines the
  ``first'' ancestors, that is those ancestors that themselves do not have
  ancestors.  Thus, the algorithm determines by itself, in a single sweep, how
  far back into the past the simulation must go.  This contrasts with usual
  CFTP schemes where algorithms may have to be iterated several times, going
  further and further into the past, until coalescence is achieved.

\end{itemize}

The relation ``being ancestor of'' induces a backwards in time
\emph{contact/oriented percolation} process.  The algorithm is applicable as
long as this oriented percolation process is subcritical.  This implies the
following limitations of our scheme:
\begin{itemize}
\item It works at \emph{low density} of objects, at least in infinite volume.
  It may work at higher densities in finite volume, but we have not pursued
  this investigation.
  
\item The birth-rate of objects must be uniformly bounded. This is necessary
  to guarantee the existence of the dominating free birth-and-death process.
\end{itemize}

Our algorithm does not rely on any type of monotonicity.  Therefore for
monotone systems our algorithm probably needs to go further back into the past
than dominated CFTP or other schemes that exploit monotonicity.  This loss in
efficiency could be compensated, at least partially, by the ``economy'' and
``single-sweep'' features mentioned in (i) and (ii) above.  As Prof.\ Kendall
pointed out (private communication), dominated CFTP constructions ``are
wasteful in that they simulate past grains without regard to whether or not
they are in the relevant percolation cluster, but efficient in that they use
some kind of monotonicity to detect whether or not one needs to investigate
further back into the past.''  As a counterpart, insensitivity to monotonicity
amounts to generality and versatility.  In particular, this versatility could
be used to offset the limitation imposed by the low-density constraint.

As an illustrative analogy, let us present a parallel with what happens in
statistical mechanics, where studies usually rely on two types of methods: (I)
those based on correlation inequalities, and (II) those expansion based.
(These are not the only methods, others include exact solutions and more
abstract arguments based on compactness or convexity.)  Types (I) and (II) are
mutually complementary.  Correlation inequalities yield very strong results,
often valid over whole regions in parameter space.  Nevertheless, these
results are rigid in that their validity depends on very precise symmetry
(monotonicity) properties that can be easily destroyed by even infinitesimal
perturbations.  In contrast, expansion-based techniques are very versatile and
robust.  While their {\it a priori} range of convergence is limited ---it is
restricted to low densities or high temperatures---, it is often the case that
suitable changes of variables place other regions of parameter space within
the scope of expansion methods.  For instance, expansion-based studies of
high-temperature spin models work with (interacting) random walks (see, e.g.\ 
Dobrushin, 1996).  Closer to, and above, the critical temperature alternative
expansions are available (Olivieri and Picco, 1990; Fern\'andez, Fr\"ohlich
and Sokal, 1992).  At low temperature the right variables are the contours,
and to get closer to (but below) the critical temperature coarse grained
contours are needed (Gaw\c edzki, Koteck{\`y} and Kupiainen, 1987).  In fact, the
belief is that there always are ``good'' variables that make everything
diluted enough.

Existing dominated-CFTP algorithms are reminiscent of methods of type (I).
They are very effective and apply for large intervals of rates.  But they are
also very specific, small alterations in the models could affect monotonicity
and render an algorithm inapplicable.  Our scheme could, perhaps, play a role
similar to the studies of type (II) for simulation purposes.  In this sense,
it is crucial that monotonicity requirements be absent from the procedure.
The goal is to change variables so to fall into a low-rate Poissonian
(birth-and-death) process.  Such changes will in general destroy any (obvious)
monotonicity property.  As an example, our scheme is capable to deal with
Peierls contours, and hence to provide an exact-sampling algorithm for the
low-temperature Ising model (of course, it is ``exact'' modulo time and space
user-impatience).  This is a region inaccessible to pre-existing algorithms.
The use of other random objects (see the end of the previous paragraph) could
yield analogous algorithms for other regions of the phase diagram.

The comparison of our algorithm with expansion methods is, in fact, more than
just an analogy.  Its theoretical basis has been used to construct an
alternative to usual expansion methods in statistical mechanical (Fern\'andez,
Ferrari and Garcia, 1998 and 2001).  This alternative has a provable region of
validity that exceeds that of usual cluster-expansion treatments.

For the sake of completeness we start with the definitions of the most
conspicuous space processes whose distributions we can perfect-simulate
(Section \ref{sec:space}). Examples include area- and perimeter-interacting
point process, invariant measures of loss networks, the random cluster model
and the contour representation of the ``$+$'' or ``$-$'' Ising measures at low
temperature.  Its relation with birth-and-death processes is discussed
immediately after (Section \ref{sec:life}), together with the basic simulation
approach for the latter.  The perfect simulation scheme is finally presented
in Section \ref{sec:perfect}.  Its central piece is the time-backward
construction of the clan of ancestors of a Poissonian birth-and-death process.

\section{Point processes} \label{sec:space}

Let $\G$ be a measurable space and $\nu$ a Radon measure on $\G$. Typically
$\G$ is $\R^d$, $\Z^d$, $\R^d\times \G'$ or $\Z^d\times \G'$, where $\G'$ is a
set of ``animals'' or ``marks''. Let
\[
{\mathcal S}=\{\xi\in\N^\G\,:\, \xi(\ga)> 0 \hbox{ only for a countable set of
}\ga\in\G\}
\]
A point process is a random element $N\in\mathcal S$. We denote with $\mu$ the
law of a point process $N$. $N$ is interpreted either as a random
configuration of points or a random counting measure on $\G$.

\paragraph*{Poisson Process} 
The first example is a Poisson process on $\G$ with intensity measure
$\nu$. Its law is characterized by
\[
\mu^0(N: N(B) = k) = {e^{-\nu(B)}\nu(B)^k/k!}
\]
for measurable $B\subset \G$; besides, under $\mu^0$ $N(B_i)$ are independent
if $B_i$ are disjoint.  When $\G=\R^d$ and $\Lambda\subset\R^d$ we call
$\mu^0_\Lambda$ the law of $N^0\cap\Lambda$. We call a Poisson process on $\G=
\R^d$ homogeneous when $\nu(\Lambda)$ is a function of $\ell(\Lambda)$, the
Lebesgue measure of $\Lambda$. Similarly, when $\G= \Z^d$, the process is
called homogeneous when $\nu(\Lambda)$ is a function of $|\Lambda|$, the
number of points in $\Z^d\cap\Lambda$. In this case, the intensity is
proportional to the Lebesgue (respectively, counting) measure and the factor
of proportionality is called the \emph{rate} which equals $\nu(\Lambda)$ for
any $\Lambda$ with unit Lebesgue measure (resp. counting measure).

\noindent{\bf Finite total rate.} For future purposes we consider the
case $\nu(\R^d\times \R^+)<\infty$; we interpret the last coordinate as time.
One can compute the distribution of the (not necessarily finite) time
$\tau_1$, the smallest time-coordinate of the points (if any) of the process.
Indeed, calling $N$ the point Poisson process with rate $\nu$, for $0\le t \le
\infty$,
\begin{equation}
  \label{t10}
  \P(\tau_1>t) \;=\; \P(N(\R^d\times[0,t))=0)
\;=\; \exp(-\nu(\R^d\times[0,t))\;.
\end{equation}
In the case of one-dimensional processes ($d=0$) the above reads
\begin{equation}
  \P(\tau_1>t) \;=\; \P(N([0,t))=0) 
\;=\; \exp(-\nu[0,t))\label{t11}\;.
\end{equation}

In this paper we consider only point processes that are absolutely continuous
with respect to a Poisson process with law $\mu^0$. The law of these processes
is characterized by
\[
\mu(dN) = \Psi(N) \mu^0(dN)
\]
where $\Psi$ is the Radon-Nikodim derivative of $\mu$ with respect to $\mu^0$. 

A Poisson process that appears in the literature is the \emph{germ-grain}
Poisson process. In this case $\G = \R^d\times\mathcal B^0(\R^d)$, where
$\mathcal B^0(\R^d)$ is the set of compact Borel sets of $\R^d$. For each
alive germ $x\in\R^d$, $g\in\mathcal B^0(\R^d)$ is the associated grain.
Assume that the grains are determined by a random variable independent of the
rest, given by a certain probability distribution $\pi_x$, which may depend on
the germ location $x$.  The intensity $\nu$ is defined by
\begin{equation}
  \label{l40}
  \nu(d(x,g)) \;=\; f(x)\,\pi_x(dg)\,dx\;.
\end{equation}
where $f(x)$ is the intensity of germs.

\paragraph*{Area-interaction point processes} 
\label{ss.area}

These processes have been introduced by Baddeley and Van Lieshout (1995). This
is a germ-grain process as defined above, but the grain shape is fixed and
equal to a compact convex $G\subset\R^d$. We only need to keep track of the
germs, so $\G=\R^d$. The intensity $\nu$ is defined by $\nu(dx) = \kappa dx$,
$\kappa$ is a positive real number.  The intersections of the grains determine
a weight that corrects the otherwise Poissonian distribution of germs. The
process is absolutely continuous with respect to the Poisson process $\mu^0$
with intensity $\nu$.  The law of the area-interaction process is defined for
bounded windows $\Lambda\subset\R^d$ by
\begin{equation}
\label{eq:300}
\mu_\Lambda (dN) \;=\; {\phi^{-m_d(N \oplus G)} \over
 Z_\Lambda(\kappa,\phi) }\,\mu^0_\Lambda(dN) \;,
\end{equation}
where $\mu^0_\Lambda$ is the law of the unit Poisson process in the box
$\Lambda$, $\phi$ is a positive parameter, $Z_\Lambda(\phi)$ is a normalizing
constant and $N \oplus G$ is the \emph{coverage process} given by
\begin{equation}
\label{eq:p4}
N \oplus G \;:=\; \bigcup_{x \in N}\, \{x + G\}\;.
\end{equation}

\paragraph*{Strauss process}
The setup is the same as the area interaction process, but now the unit
Poisson process is weighted according to an exponential of the number of pairs
of points closer than a fixed threshold $r$. The measure is defined by
\begin{equation}
\label{eq:300ca}
 \mu_\Lambda(dN) = \frac{1}{Z_\Lambda} e^{\beta_1 N(\Lambda) + \beta_2
 S(N,\Lambda)}\, \mu^0_\Lambda(dN) 
\end{equation}
where $S(N,\Lambda)$ is the number of unordered pairs such that $\|x_i-x_j\| <
r$. The case $\beta_2>0$ was introduced by Strauss (1975) to model the
clustering of Californian red wood seedling around older stumps. However,
\reff{eq:300ca} is not integrable in that case (see Kelly and Ripley (1976)).

\paragraph*{Low-temperature Ising model}\label{ss.ising}
The well-known \emph{Peierls contours} allow to map the ``+'' or ``$-$''
measures of the ferromagnetic Ising model at low temperature into an ensemble
of objects ---the contours--- interacting only by perimeter-exclusion.  See,
for instance, Section 5B of Dobrushin, 1996, for a concise and rigorous
account of this mapping.  The (discrete) set $\G$ consists of contours; these
are hypersurfaces formed by a finite number of $(d-1)$-dimensional unit cubes
---\emph{links} for $d=2$, \emph{plaquettes} for higher dimensions--- centered
at points of $\Z^{d}$ and perpendicular to the edges of the dual lattice
$\Z^d+({1\over 2},\cdots,{1\over 2})$. To each contour one can assign an
``origin'' in $\Z^d$ and say that two contours are equivalent if they coincide
after a translation of the origin. Calling $\G'$ the set of contours modulus
this class of equivalence, the set $\G$ can be expressed by
$\G=\Z^d\times\G'$, where the first coordinate corresponds to the origin and
the second to the ``shape'' of the contour.  Call two plaquettes
\emph{adjacent} if they share a $(d-2)$-dimensional face.  A set of
plaquettes, $\gamma$, is \emph{connected} if for any two plaquettes in
$\gamma$ there exists a sequence of adjacent plaquettes in $\gamma$ joining
them.  The set $\gamma$ is \emph{closed} if every $(d-2)$-dimensional face is
covered by an even number of plaquettes in $\gamma$.  \emph{Contours} are
connected and closed sets of plaquettes.  For example, in two dimensions
contours are closed polygonals.  Two contours $\ga$ and $\theta$ are said to
be \emph{compatible} if no plaquette of $\ga$ is adjacent to a plaquette of
$\theta$.  In two dimensions, therefore, contours are compatible if and only
if they do not share the endpoint of a link.  In three dimensions two
compatible contours can share vertices, but not sides of plaquettes.  Ising
spin configurations in a bounded region with ``+'' (or ``$-$'') boundary
condition are in one-to-one correspondence with families of pairwise
compatible contours.

Let the \emph{compatibility matrix} $I:\G\times\G\to \{0,1\}$ be defined by
\begin{equation}
  \label{eq:c3}
  I(\ga,\theta) \;=\; \left\{
\begin{array}{ll} 
0, & \hbox{if $\ga$ and $\theta$ are compatible}\\
1, & \hbox{otherwise}
\end{array}\right.
\end{equation}
The ``Poisson process'' $\mu^0_\Lambda$ in $\N^\G$ is the product of Poisson
measures whose $\gamma$-marginal is Poisson with mean
\[w(\ga) := \exp(-\beta|\gamma|)\]
for $\gamma\subset\Lambda$. Here $|\gamma|$ stands for the number of
plaquettes of $\gamma$. The intensity measure $\nu$ is discrete: $\int_B
d\nu=\nu(B) = \sum_{\ga\in B} w(\ga)$.

Let $\mu_\Lambda$ be the measure defined by: for $\xi\in\{0,1\}^\G$, such that
$\xi(\ga)\le \one\{\ga\subset\Lambda\}$, 
\begin{equation}
  \label{c10}
  \mu_\Lambda(\xi) \,= {1\over Z_\Lambda}\,
  \Bigl(\prod_{\ga,\th:\xi(\ga)\xi(\th)=1} 
  [1-I(\ga,\th)]\Bigr)\,\mu^0_\Lambda(\xi) 
\end{equation}
where $\beta$ is a positive parameter called inverse temperature.  The factor
$Z_\Lambda$ is just the normalization.

\paragraph*{The random cluster model}\label{ss.rc}

Consider $\Lambda\subset\Z^d$ and let $\BB(\Lambda) := \{(x,y)\in
\Lambda\times\Lambda: |x-y|=1\}$ the set of bonds of $\Lambda$. A bond
configuration $\zeta\in\{0,1\}^{\BB(\Lambda)}$ is a function from
$\BB(\Lambda)\to\{0,1\}$. Bonds assigned 1 are called \emph{open}, otherwise
\emph{closed}. A \emph{cluster} of $\zeta$ is a set of sites connected with
open bonds; sites surrounded only by closed bonds are clusters of size 1. Let
$p\in[0,1]$ and $q>0$ be parameters and define the finite volume measure
\begin{equation}
  \label{rcm}
  \varphi_\Lambda(\ze) = {1\over Z_\Lambda(p,q)} \,
  p^{O(\zeta)}\,(1-p)^{C(\zeta)}\, q^{L(\ze)} 
\end{equation}
where $O(\ze)$ is the number of open bonds of $\ze$ in $\Lambda$, $C(\zeta)$
is the number of closed bonds and $L(\ze)$ is the number of clusters. The
constant $Z_\Lambda(p,q)$ is the normalization. In other words,
$\varphi_\Lambda$ is absolutely continuous with respect to the product measure
on $\BB(\Lambda)$ with parameter $p$, with Radon-Nikodim derivative
$q^{L(\ze)}/Z_\Lambda(p,q)$. This model has been introduced by Fortuin and
Kasteleyn (1972); a review can be found in Grimmett (1995).

Taking the connected sets of bonds as the basic objects, this model can be
written as in \reff{c10}. More precisely, if one says that two sets of bonds
are \emph{incompatible} whenever they share some vertices, and takes $\G=\{
\ga\subset\BB(\Lambda)\,:\, \ga$ is finite and connected$\}$, the probability
weights of the model can be written in the form \reff{c10} with the Poisson
means
\begin{equation}
  \label{wei}
  w(\ga) \;=\;
 \Bigl({p\over 1-p}\Bigr)^{B(\ga)}\,\Bigl({1\over q}\Bigr)^{V(\ga)-1}\;.
\end{equation}
Here $B(\ga)$ is the number of bonds of $\ga$ and $V(\ga)$ the number of
vertices in the extremes of the bonds of $\ga$. [That is, $V(\ga)=\#\{x\in
\Lambda \hbox{ such that } (x,y)\in\ga \hbox{ or } (y,x)\in\ga$ for some
$y\in\ga\}$.]  The transformation $ Y \,:\, \{0,1\}^{\BB(\Lambda)} \,
\longrightarrow \, \{0,1\}^{\G}$ defined by
\[
Y(\zeta)(\gamma)=1 \quad \Longleftrightarrow \quad
\hbox{$\ga$ is a maximally connected set of open bonds of $\zeta$}
  \label{eq:wei.1}
\]
satisfies
\begin{equation}
  \label{mf}
  \varphi_\Lambda(\zeta) \,=\, \mu_\Lambda(Y(\zeta))\,.
\end{equation}

\section{Birth-and-death processes}
\label{sec:life}

\subsection{Definition and examples}
The common feature linking all the spatial processes described in the previous
section is that \emph{all these distributions can be realized as invariant
  measures of spatial interacting birth-and-death processes}.

We consider the state space ${\mathcal S}$ of point configurations on $\G$
with a Radon measure $\nu$ as in Section 2.

The \emph{free} birth death process is characterized by the fact that
individuals are born at intensity $\nu$ and last for a random time
exponentially distributed of mean one.  The generator of the free process is
\begin{eqnarray}
  \label{gen0}
  A^0F(\eta)&=& \int_\G \nu(d\ga) \,
  [F(\eta+\delta_\ga)- F(\eta)]\nonumber\\
&&\qquad+\; \sum_{\ga\in\G:\eta(\ga)>0}\, \eta(\ga)\,
  [F(\eta-\delta_\ga)- F(\eta)]  
\end{eqnarray}
Here $\delta_\ga$ is the configuration with only one point at $\ga$ and
$(\eta+\xi)(\th) = \eta(\th)+\xi(\th)$ (coordinatewise sum). The invariant
(and reversible) measure for the free process is the Poisson process $\mu^0$
with intensity $\nu$.

Let $\mu$ be a measure absolutely continuous with respect to $\mu^0$ with
Radon-Nikodim derivative $\Psi$: $d\mu(\eta) = \Psi(\eta) d\mu^0(\eta)$. Define
\begin{eqnarray}
  \label{gen}
  AF(\eta)&=& \int_\G \nu(d\ga) \,{\Psi(\eta+\delta_\ga)\over \Psi(\eta)}\,
  [F(\eta+\delta_\ga)- F(\eta)]\nonumber\\
&&\qquad+\; \sum_{\ga\in\G:\eta(\ga)>0}\, \eta(\ga)\,
  [F(\eta-\delta_\ga)- F(\eta)]  
\end{eqnarray}
The difference with the free process is that in the interacting process the
rate of birth $\nu(d\ga)$ is corrected with the quotient
$\Psi(\eta+\delta_\ga)/ \Psi(\eta)$, while the rate of death remains
unaltered. The measure $\mu$ is reversible for the process with generator $A$.
To better interpret this dynamics assume
\[
\Delta_\Psi \;:= \;\sup_{\eta,\ga} {\Psi(\eta+\delta_\ga)\over
  \Psi(\eta)}\;<\;\infty
\]
and define $M: \G\times\mathcal S \to[0,1]$ by 
\begin{equation}
  \label{mmm1}
  M(\ga|\xi)= {\Psi(\xi+\delta_\ga)\over \Delta_\Psi\,\Psi(\xi)} 
\end{equation}
If $\bar\nu = \Delta_\Psi\nu$, the generator can be rewritten as 
\begin{eqnarray}
  \label{genp}
  AF(\eta)&=& \int_\G \bar\nu(d\ga) \,M(\ga\vert \eta)\,
  [F(\eta+\delta_\ga)- F(\eta)]\nonumber\\
&&\qquad+\; \sum_{\ga\in\G:\eta(\ga)>0}\, \eta(\ga)\,
  [F(\eta-\delta_\ga)- F(\eta)]  
\end{eqnarray}
This dynamics has the following interpretation. When the current configuration
of objects is $\xi$, object $\gamma$ attempts to be born with rate
$\bar\nu(d\ga)$ and is effectively born with probability $M(\ga|\xi)$. The
death rate of any object is one.

The interaction $M$ induces naturally the notion of \emph{incompatibility}
between individuals. This is a not necessarily symmetric matrix
$I:\G\times\G\to\{0,1\}$ defined by
\begin{equation}
  \label{com1}
  I(\ga,\th) := \one\Bigl\{\sup_{\grains}\Bigl\{\Bigl|M(\ga\,\vert\,\grains)-
  M(\ga\,\vert\,\grains+\delta_{\th})\Bigr|\Bigr\}>0\Bigr\}\;,
\end{equation}
where $\delta_\th$ is the configuration having unique individual $\th$ and the
supremum is taken over the set of those $\xi$ such that $\xi$ and
$\xi+\delta_\th$ are in the set of configurations (either $\{0,1\}^\G$ or
$\N^\G$).  The function $I(\ga,\th)$ indicates which individuals $\th$ may
have an influence in the birth-rate of the individual $\ga$. In the case of
the spatial point processes described above the matrix is symmetric and given
by $I((x,g),(x',g'))=\one\{(x+g)\cap (x'+g')\neq\emptyset\}$. If
$I(\ga,\th)=1$, ---that is, if the presence/absence of $\th$ modifies the rate
of birth of $\ga$--- we say that $\th$ is \emph{incompatible} with $\ga$. For
the Ising and random-cluster models one recovers \reff{eq:c3}.

\paragraph*{Area interaction Point processes}
In the repulsive ($\phi<1$) point process \reff{eq:300} we have $\Delta_\Psi=
\kappa\phi^{-m_d(G)}$ and 
\begin{equation}
  \label{eq:l42}
  \bar\nu(dx) \;=\; \kappa\,\phi^{-m_d(G)}dx\;,
\end{equation}
\begin{equation}
  \label{eq:l43}
  M(x\vert\grains) = \phi^{m_d(G) - m_d((x+ G) \setminus (\xi
  \oplus G))} 
\end{equation}
For the attractive ($\phi>1$) case, 
\begin{equation}
  \label{eq:l44}
 \bar\nu(dx) \;=\; \kappa\,dx
\end{equation}
\begin{equation}
  \label{eq:l43.1}
  M(x\vert \xi) = \phi^{- m_d((x + G)\setminus (\grains \oplus G))}\;.
\end{equation}

\paragraph*{Continuous unbounded one-dimensional loss network}
\label{ss.one}

A loss network models, for instance, the occurrence of calls in a
communication network.  Kelly (1991) reviews several discrete regimes and
introduces the following continuous generalization. Callers are arranged along
an infinitely long cable and each call between two points $s_{1}$, $s_{2} \in
\R$ on the cable involves just the segment between them.  The cable has the
capacity to carry simultaneously up to $C$ calls past any point along its
length.  Hence, a call attempt between $s_{1}$ and $s_{2} \in \R$, $s_{1} <
s_{2}$, is lost if past any point of the interval $[s_{1}, s_{2}]$ the cable
is already carrying $C$ calls. Calls are attempted with initial (leftmost)
point following a space-time Poisson process with intensity $f(x) dx$, and
(space) lengths given by a distribution $\pi$, independent of its leftmost
point, with finite mean $\rho_1$.  The holding time of a call has exponential
distribution with mean one.  The location of a call, its length and its
duration are independent.

In this case, the germs ($x$) are the leftmost points of calls and the grains
($g$) are segments with random lengths. This process can be viewed as a
spatial birth and death process where we can take
\begin{equation}
  \label{eq:l42a}
  \bar\nu(d(x,g)) \;=\; f(x) \,dx \,\pi(dg) \;,
\end{equation}
and denoting $\xi(u):=$ number of calls using point $u$, 
\begin{equation}
  \label{eq:l43a}
  M((x,g)\vert\grains) = \one{\{(\grains + \delta_{\{x+g\}})(u) \le C,
  \forall u \in \R\}}.
\end{equation}

\paragraph*{Discrete processes}

Free birth-and-death processes with a countable family of individuals $\G$ are
simply the product of independent birth death processes labeled by each
$\gamma\in\G$, with birth rates $w(\gamma)$ and death rate equal to the number
of alive individuals.  Such a process exists without any requirement on the
weights $w(\gamma)$; it is ergodic and its invariant distribution is the
product of Poisson laws with mean $w(\gamma)$.

The (discrete) loss networks, the contour model and the animal version of the
random cluster model of Section 2 are processes of this form where, in fact,
the matrix $M$ takes only two values, 0 and 1.  That is, the interaction
imposes a deterministic constraint.  In particular, the interaction terms of
the Ising-model, random cluster and loss networks have a simple product form
\begin{equation}
  \label{eq:l32}
  M(\ga\vert\xi) \;=\;  \prod_{\theta:\xi(\theta)\neq 0} [1-I(\ga,\theta)]\;.
\end{equation}
Indeed, the Radon-Nikodim derivative is one for allowed configurations, hence
the denominator in \reff{mmm1} is one.

\subsection{Graphical construction}

We proceed to the construction of the probability space where both the free
process and the interacting birth death process will be constructed. Consider
the countable family of random quartets
$\{(\Gamma_i,T_i,S_i,Z_i)\,:\,i\in\J\}$, with $\Gamma_i\in\G$, $T_i,S_i\in\R$,
$Z_i\in[0,1]$ such that:

\begin{itemize}
  
\item The process $\bigl\{(\Gamma_i,T_i) \,,\,i\in\J \bigr\}$ is a Poisson
  process on $\G\times\R$ with mean measure $\nu\times\ell$; $\ell$ is the
  Lebesgue measure in $\R$.  This process determines the times and type of
  attempted births of individuals.
 
\item $S_i$ is exponentially distributed with mean 1. This variable will
  determine the lifetime of the $i$th attempted birth.
  
\item $Z_i$ is uniformly distributed in $(0,1)$. This variable is called the
  \emph{flag} or \emph{mark} of the $i$th attempted birth and will be used
  together with the function $M(\cdot|\cdot)$ to decide if the attempted birth
  is actually a real birth.
\end{itemize}

Each triplet $(\Gamma_i,T_i,S_i)$ can be visualized as a
\emph{cylinder} of (space) \emph{basis} $\Gamma_i$, birth time $T_i$
and lifetime $S_i$.  The random set of marked cylinders is called
\begin{equation}
\label{ca:3}
\C = \Bigl\{ \Bigl(\Gamma_i \times
[T_i, T_i + S_i]\,,\,Z_i\Bigr) \;,\; i\in\J \Bigr\}\;.
\end{equation}
For a generic marked cylinder $C = (\Gamma\times[t,t+s],z)\in\C$, denote
$\birth(C) = t$, $\death(C)=t+s$, $\life(C)= [t,t+s]$, $\basis(C)=\Gamma$ and
$\flag(C)=z$.

\paragraph*{The free process} 
The construction of $\C$ is time-translation invariant. Call
\begin{equation}
\label{eq:l45}
\xi_t(\ga):= \#\Bigl\{C\in \C\,:\, \basis(C) =
\ga\,;\, \life(C)\ni t\Bigr\}
\end{equation}
the set of individuals forming the sections of $\C$ at time $t$. All attempted
births are actual births in this case. $\xi_t(\ga)$ will be at most 1 in the
continuous case, but could be bigger in the discrete case. $(\xi_t:t\in\R)$
constitutes a stationary free birth-death process with generator \reff{gen0}.
The marginal law of $\xi_t$ is $\mu^0$, the Poisson process with intensity
$\nu$.

Likewise, one can define the free process on $\R^d\times [0,\infty)$ with
initial configuration of individuals
\begin{equation}
\label{ca:4}
\grains_0 \, := \, \{\gamma^0_1, \gamma^0_2, \ldots\}\;.
\end{equation}
For this associate cylinders to the initial configuration:
\begin{equation}
\label{ca:5}
\C_0(\grains_0) := \Bigl\{\bigl( \gamma^0_i  \times
[0, 0 + S^0_i]\,,\,Z^0_i\bigr) \;,\, i\in\Z \Bigr\}
\end{equation}
where $S^0_i$ and $Z_i^0$ are independent and independent-of-everything random
variables whose distributions are, respectively, Exp$(1)$ and $U(0,1)$.
Define the subset of cylinders born between $0$ and $t$:
\[
\C_{[0,t]} := \{C\in \C\,:\, \birth(C)\in[0,t]\}
\]
Then, the process defined at time $t$ by
\begin{equation}
\label{l450}
\xi_t(\ga):= \#\Bigl\{C\in \C_{[0,t]}\cup\C_0(\grains_0)\,:\, \basis(C) =
\ga\,;\, \life(C)\ni t\Bigr\} 
\end{equation}
has initial configuration $\xi_0$ and generator \reff{gen0}.

\paragraph*{Interacting processes} 
\label{sub:sbdp}

The absolute continuity with respect to the free process, embodied in the
generator \reff{genp}, suggests a simple alteration to the previous
construction to pass to an interacting birth-and-death process: The
\emph{attempted births} become actual births only if an additional (generally
stochastic) test is passed.  This test is determined by the factor $M$ of the
rate densities.  The interacting process is, therefore, obtained as a
``thinning'' or ``trimming'' of the free process.

The formalization of this intuitively simple idea is easy for finite windows,
but more delicate for the infinite-volume process.  We discuss the former case
first.

\subsubsection{Finite-volume construction}
\label{ss.finite}

To construct a birth-and-death process $\xi_t$ with rate density
$\bar\nu(d\gamma)\,M(\gamma\vert\grains)$, for individuals within a finite
space-region $\Lambda$ and for a finite time interval $[t_0,t_{\rm fin}]$, one
proceeds as follows:

\begin{enumerate}
  
\item Run the free process with rate density $\bar\nu$ starting from the
  initial cylinders $\C_0$.  If $M(\gamma\vert\grains)$ is deterministic
  ---for instance forbidding individuals to overlap--- the initial
  configuration is assumed to satisfy the corresponding constraint.
  
\item Each death happening before reaching an event of the free process causes
  the corresponding updating of $\grains_t$, by taking the corresponding
  individual out of $\xi_t$.
  
\item When the free process yields a first event $(\gamma_1,t_1,s_1,z_1)$,
  this event is considered an attempted birth.  To decide, one looks to the
  set $\grains_{t_1-}$ of alive individuals ($\grains_{t_1-}$ is equal to
  $\grains^0$ minus the initial individuals with lifetime smaller than $t_1$).
  If
\begin{equation}
\label{eq:t50}
z_1 \;<\; M(\gamma_1\vert\grains_{{t_1}-})
\end{equation}
the cylinder is allowed to be born and the individual $\gamma_1$ is included
in the configuration $\grains_{t_1}$; otherwise it is ignored and
$\grains_{t_1}$ is set equal to $\grains_{t_1-}$.

\item Now iterate the procedure, that is, repeat the previous two steps
  shifting subscripts $1\to 2$ and $0\to1$.  Continue in this way until
  reaching an attempted birth beyond $t_{\rm fin}$.

\end{enumerate}

\subsubsection{Two-sweep finite-volume construction}
\label{ts.two}

The visualization in terms of cylinders suggests an alternative implementation
as a two-sweep scheme: In the first sweep one generates free cylinders by
running the free process from $t_0$ to $t_{\rm fin}$, while in the second
sweep a decision is made on which cylinders are kept and which are erased.
The set of \emph{kept cylinders} includes, by definition, all initial
cylinders while successive additions must pass the test \reff{eq:t50}.  We
call $\K_{[0,t]}(\Lambda,\xi_0)$ the resultant set of kept cylinders in the
construction of Section \ref{ss.finite}. The configuration of the process at
time $t$ with initial configuration $\xi_0$ is then given by the projection of
the bases of the alive kept cylinders at that time:
\begin{equation}
  \label{449}
  \eta_t(\gamma) =  
  \#\Bigl\{C\in\K_{[0,t]}(\Lambda,\xi_0)
  \hbox{ with basis } \gamma \hbox{ and alive at } t \Bigr\} \;.
\end{equation}

\subsubsection{Finite-volume time-stationary construction}
\label{ts.finite}

The construction can be also performed in a stationary manner for $t\in \R$.
Indeed, since in a finite window the number of alive individuals is finite
(with probability one), there exist random times $\{\tau_j\in\R: j\in\Z\}$,
such that (a) $\tau_j \to \pm\infty$ for $j\to\pm\infty$ and (b) $\bigl(\cup_i
[T_i,T_i+S_i]\bigr) \cap \bigl(\cup_j\{\tau_j\}\bigr) = \emptyset$.  In words,
at each $\tau_j$ no cylinder is alive. The above selection of kept cylinders
can then be performed independently in each of the random intervals
$[\tau_i,\tau_{i+1})$.  This stationary construction is particularly useful to
study properties of the invariant measure $\mu_\Lambda$. In fact, calling
$\K(\Lambda)$ the (time stationary random) set of kept cylinders, the law of
\begin{equation}
  \label{447}
  \eta_t(\gamma) =  
  \#\Bigl\{C\in\K(\Lambda)
  \hbox{ with basis } \gamma \hbox{ and alive at } t \Bigr\}
\end{equation}
is exactly $\mu_\Lambda$. 

\subsubsection{Infinite-volume construction}
\label{ss.infinite}

None of the finite-volume procedures discussed above can be directly
implemented to construct the process in infinite volume.  On the one hand, the
scheme proposed for finite time intervals is not applicable to infinite volume
because it is not possible to decide which is the first mark in time.  On the
other hand, the stationary construction is also not feasible because in
infinite volume there are cylinders alive at all times.  This last objection,
however, may play no role if one only focuses on a family of cylinders
intersecting a \emph{finite} set and tries to decide which of them should be
erased and which ones kept.

According to the previous discussion, to decide whether a cylinder $C\in \C$
is kept, one has to look at the set of cylinders $C'$ (born before $C$ and)
alive at the birth-time of $C$ whose basis are incompatible with the basis of
$C$ in the sense of \reff{com1}. Let us call this set the first generation of
\emph{ancestors} of $C$ and denote it $\A^C_1$. Once we determine which of
these ancestors are alive, the decision on whether to keep $C$ or not requires
only a single application of the test \reff{eq:t50}.  However, to decide which
of these ancestors are alive we have to work with the second generation of
ancestors of $C$, that is, with the ancestors of the ancestors.  Recursively,
we find ourselves having to deal with all generations of ancestors of $C$.
Let us call the union of all generations of ancestors of $C$ the \emph{clan of
  ancestors} of $C$, and denote it $\A^C=\cup_{n\ge 1}\A^C_n$, where $\A^C_n$
is the set of ancestors in the $n$th generation.  These sets may contain
cylinders in $\C_{\zero}(\xi_0)$. The procedure for deciding whether to keep
or to erase $C$ can be univocally defined if \emph{the clan of ancestors of
  $C$ is finite}.

This picture makes it apparent that an infinite-volume process given initial
starting conditions (i.e. for a \emph{finite time-interval}) exists as long as
there are no explosions, that is, as long as no cylinder can develop
infinitely many ancestors in a finite time.  Furthermore, there exists a
unique stationary process (for \emph{infinite time-intervals}) if \emph{all}
clans of ancestors are finite with (free-process) probability one.

\begin{thm}\label{thm1}\mbox{}

\begin{itemize}
\item[(i)] If with probability one $\A^C\cap\C_{[0,t]}$ is finite for every
  cylinder $C$ alive at time $t$, for any $t>0$, then the birth-and-death
  process with the generator \reff{gen} and initial condition $\eta_0$ is
  obtained by performing the two-sweep construction of Section \reff{ts.two}
  on each set $(\A^C\cap\C_{[0,t]})\cup\C_0(\eta_0)$ and taking the
  projections
\begin{equation}
  \label{4440}
  \eta_t(\gamma) =  
  \#\Bigl\{C\in\K_{[0,t]}(\eta_0) \hbox{ with
  basis } \gamma \hbox{ and alive at } t \Bigr\}\;.
\end{equation}
where $\K_{[0,t]}(\eta_0)$ is the resulting set of kept cylinders.
\item[(ii)] If with probability one $\A^C$ is finite for every cylinder $C$,
  then the stationary birth-and-death process with the generator \reff{gen}
  can be constructed for $t\in\R$ by performing the two-sweep construction of
  Section \reff{ts.two} on each set $\A^C$ and taking the projections
\begin{equation}
  \label{4441}
  \eta_t(\gamma) =  
  \#\Bigl\{ C\in\K \hbox{ with
  basis } \gamma \hbox{ and alive at } t \Bigr\}\;.
\end{equation}
Moreover, the marginal distribution of $\eta_t$ is the stationary measure
$\mu$.
\end{itemize}  
\end{thm}

A proof of this theorem is presented in the Appendix.

\section{Oriented percolation and branching processes}
\label{pbp}
To determine the conditions allowing the construction of Theorem \ref{thm1},
we point out that the relation ``being ancestor of'' gives rise to a model of
\emph{oriented-percolation}.  We call it \emph{backwards oriented percolation}
to emphasize the fact that it is defined by only looking into the past.  The
finite-time construction is possible if there is no cluster with infinitely
many members in a finite time slice, while the feasibility of the
infinite-time construction requires the absence of a percolation cluster
reaching to time $-\infty$.  As usual in oriented percolation problems, it is
useful to work with a majorizing \emph{multitype branching process}. In this
process the offspring distribution of a cylinder $C$ has the same (marginal)
law as the distribution of $\A^C_1$, but the branches behave independently.
The problem is then reduced to determine conditions guaranteeing the
finiteness of the clan of branching ancestors. Hence it is sufficient to show
in the finite-time case that the branching process does not explode, while in
the infinite-time case we need to prove that the branching process is
sub-critical. Factorization makes these tasks easier.

Let us give sufficient conditions on the dominating branching for the
different processes listed in Section \ref{sec:space}.

\paragraph*{Discrete processes} 
\label{sub:gcni}
If the family of individuals $\G$ is countable, the free birth-and-death
process is the product over $\gamma\in\G$ of independent marked Poisson
processes.  The construction of the interacting processes is an obvious
adaptation of the procedure of the continuous case. For the infinite-volume
process one relies on the properties of the backwards oriented percolation
model of cylinders defined by the oriented bonds $C\to C'$ if $C'$ is an
ancestor of $C$, that is if the basis of $C$ and $C'$ intersect and $C'$ is
alive when $C$ is born. Let $m(\ga,\th)$ be the mean number of cylinders of
basis $\theta$ in the first generation of a cylinder of basis $\gamma$. These
are cylinders born at negative times $-t$ and have a lifetime at least $t$, so
they survive to intersect the grain born at time zero.  Its average number is,
therefore,
\begin{equation}
  m(\ga,\th) = w(\th) I(\ga,\th)\,
\Bigl[\int_{-\infty}^0 dt\, \int_{t}^{\infty} ds\, \ee^{-s}\Bigr]
 \;=\; w(\th) I(\ga,\th)\, \;\cdot \; 1\;.
  \label{t53d}
\end{equation}
Define $m^n(\ga,\th)$ as the mean number of cylinders of basis $\th$
incompatible with a cylinder of basis $\ga$ in the $n$-th generation of
ancestors, $m^n$ is the matrix-product of $m$ by itself $n$ times.  The
condition for absence of oriented percolation is
\begin{equation}
  \label{mn1}
  \sum_{n\ge 1} \sum_\th m^n(\ga,\th) <\infty
\end{equation}
for all $\ga$. For any function ``size'' $q:\G\to\R^+$, such that $\inf_\ga
q(\ga)\ge1$, as in Lemma 5.15 of Fern\'andez, Ferrari and Garcia (2001),
calling
 \begin{equation}
  \label{mn2}
  \alpha_q\,:=\,\sup_\ga {1\over q(\ga)}\sum_\th q(\th)\, m(\ga,\th),
\end{equation}
we have
\begin{equation}
  \label{mn3}
\sum_\th m^n(\ga,\th) \,\le\, \alpha^{n}_{q}\,q(\ga)\,.
\end{equation}
The form of measuring this ``size'' depends on the process in question, but
usually there is an obvious prescription.  For instance, for the loss
networks, the Peierls contours and the random clusters model this measure is
just the length of the call, the number of plaquettes of the contour or the
number of points of the cluster. The (infinite-volume) birth-and-death process
---and hence the corresponding loss network--- exists for finite
time-intervals if $\alpha_q<\infty$, while $\alpha_q<1$ is a sufficient
condition to be an ergodic infinite-time process.

\paragraph*{Area interaction point process} Here $m(x,.)$ is a measure
on $\G=\R^d$; $m(x,dy)$ represents the rate at which cylinders with basis
centered at $y$ appear. Consider a germ $x_0$, which, by
space-time-invariance, can be placed at the origin and assumed to be born at
time zero.  Its ancestors are all cylinders whose bases involve germs located
in $\partial G:= \Bigl\{ x : (x+G) \cap G \neq \emptyset\Bigr\}$. Therefore,
as in \reff{t53d},
\begin{equation}
 m(x,\G) \;=\; \nu(\partial G) \;\cdot \; 1\;.
  \label{eq:t53}
\end{equation}
We conclude that the corresponding birth-and-death process exists
for finite times as long as $\nu(\partial G)<\infty$ and, if
furthermore,
\begin{equation}
  \label{eq:t54}
  \nu(\partial G)\;<\; 1
\end{equation}
then there is an stationary ergodic process, absolutely continuous
respect to the free process, having as invariant measure the
corresponding point process of Section \ref{sec:space}.  

The argument also works if $G\in\G$ is a random set chosen independently of
everything as in \reff{l40}. Recall $f(x)$ is the intensity of germs and
$\pi_x$ is the distribution of the grain centered in $x$. Let
\begin{equation}
  \label{mmm}
  m(G,dH)\; =\;\int dx f(x)\, I(G,H)\,\pi_x(dH)
\end{equation}
the rate at which individuals $H$ having an influence in the birth-rate of $G$
appear (see \reff{com1} for the definition of $I$). This implies that the mean
number of individuals in the first generation of ancestors of $G$ is
$m(G,\G)$. Let the ``matrix product'' $m^n$ be defined inductively by $m^1 =
m$ and
\begin{equation}
  \label{t54x}
  m^n(G,dH) := \int_{\G} m^{n-1}(G,dK)\,m(K,dH) \qquad n>1\;.
\end{equation}
As for $n=1$, $m^n(G,\G)$ is the mean number of individuals in the $n$-th
generation of ancestors of $G$. The ergodicity of the process is implied by
\begin{equation}
  \label{t54z}
  \sum_{n\ge 1}  m^n(G,\G) < \infty
\end{equation}
for all individual $G$.  A sufficient condition for \reff{t54z} is
\begin{equation}
  \label{t54a}
  \alpha_q:=\sup_G {1\over q(G)} \int_{\R^d} f(x)\, dx \int_{\G}\,
  \pi_x(dH)\; q(H) \;I(G,H) \;<\; 1\,.
\end{equation}
for some function $q:{\mathcal B}(\R^d)\to\R^+$ satisfying $\inf_Gq(G)\ge 1$.
Indeed, it can be proven as in Lemma 5.15 in Fern\'andez, Ferrari and Garcia
(2001)
\begin{equation}
  \label{666}
  m^n(G,\G)\;\le\; q(G) \,\alpha_q^n\,.
\end{equation}
Strictly speaking the above statements have been rigorously proven only for
the discretized version of the models.

\paragraph*{Loss networks}
The calls of the loss networks can be interpreted as germ-grains. For
instance, in the one-dimensional case, the germs are the leftmost points of
the calls and the grains are segments with random lengths.

A particular case where one can explicitly compute the sufficient condition is
the one-dimensional continuous loss networks of Section \ref{ss.one}.  Assume,
in general, that the leftmost points of calls appear with rate $f(x)$ and that
call lengths are given by a distribution $\pi$ independent of $x$.  We only
require the latter to have a finite mean $\rho_1$.  Consider a germ sitting at
the origin, that is a call stretching from the origin to the right, born at
time zero.  Its ancestors correspond to cylinders with sufficient lifetime and
with bases given by either calls starting at negative sites and passing
through the origin, or calls of arbitrary length originating within the sites
occupied by the initial call. Therefore, the $\alpha_q$ in \reff{t54a} for the
case $q(L)\equiv 1$ turns to be:
\begin{equation}
  \label{eq:t55a1}
 \alpha_q\, =\,  \sup_{L} \left(\int_{-\infty}^{0} \pi\{L > -x\} \, f(x) \, dx
  \,+\,  
  \int_0^L dx \, f(x)\right)\,.
\end{equation}
In the homogeneous case ($f(x) \equiv \kappa$) this gives the following
condition for ergodicity:
\begin{equation}
  \label{eq:t55a1}
  \kappa\, (\rho_1 + \sup_L L)\,<\,1.
\end{equation}
A simple computation shows that choosing $q(L) = \max(L,1)$ gives
$\alpha_q\,\le\,\kappa(\rho_2 + \rho_1+1)$, where $\rho_1$ and $\rho_2$ are
the first and second moment of the distribution $\pi$ respectively. This gives
the following sufficient condition for ergodicity
\begin{equation}
\label{eq:t54b}
\kappa(\rho_2 + \rho_1+1)\,<\,1\,.
\end{equation}
M\'{a}ric (2002) improved this bound to
 \begin{equation}
\label{eq:t54c}
\kappa(\sqrt{\rho_2} + \rho_1) <1.
\end{equation}
We remark that to obtain these conditions it was important to consider only
\emph{oriented} percolation. The analogous conditions obtained by considering
unoriented percolation of cylinders are far more restrictive.

\section{Perfect simulation of invariant measures of 
birth-and-death processes} \label{sec:perfect}

The main issue of this section is a construction of the set $\A^{\Lambda,0}$
formed by the cylinders with bases intersecting the space-time set
$\Lambda\times\{0\}$ (``cylinders alive at time 0'') and their clans of
ancestors.  This is a problem of simulation of cylinders generated by the
\emph{free} process.  Once these clans are perfectly simulated, it is only
necessary to apply the \emph{deterministic} ``cleaning procedure'', based on
the test \reff{eq:t50}, to obtain a perfect sample of the interacting process.
The scheme is feasible if these clans are finite with probability one, a fact
valid under conditions like \reff{eq:t54}, \reff{t54z} or $\alpha_q<1$, where
$\alpha_q$ is defined in \reff{mn2} for the discrete case and \reff{t54a} for
the area interaction process.

We propose a non-homogeneous time-backwards construction of the clan based on
a result proven in Section 4.5.1 of Fern\'andez, Ferrari and Garcia (2001).
It is shown there that the clan of ancestors of a family of cylinders can be
obtained combing \emph{back} in time and \emph{generating births} of ancestors
with an appropriate rate.  Alternatively, one could use the fact that the law
of $\C$ is time-reflexion invariant, to \emph{generate deaths} of ancestors.
This is simple in the area-interaction process with a fixed grain, but it is
not Markovian and more involved in the infinite case. This approach was
proposed by one of the referees and developed by Garcia (2000).

For concreteness, let us discuss our scheme for individuals living in $\R^d$
or $\Z^d$.  The birth-rate of a new cylinder to be added to the clan is equal
to the rate density of the free process multiplied by an exponential time
factor ensuring that the ancestor has a lifespan large enough to actually be
an ancestor. This time factor involves the time-distance to the birth of
existing cylinders, which can be expressed through the following function.
For a finite region $\Lambda$ and a finite set of cylinders $\H$, let the set
of bases of the potential ancestors of $\H$ and $\Lambda\times\{0\}$ be
\begin{eqnarray}
  \label{bpa}
  \G(\H,\Lambda) &:=&\Bigl\{\th\in\G\,:\, I(\basis(C'),\th)=1,\,\hbox
  { for some   }C'\in \H\Bigr\}\nonumber\\
&&\qquad\bigcup
\Bigl\{\th\in\G\,:\,\th\cap\Lambda \neq\emptyset\Bigr\}
\end{eqnarray}
and for a given individual $\th\in\G(\H,\Lambda)$,
\begin{equation}
  \label{coo}
  \ti(\H,\Lambda,\th) \;= \; 
\min\Bigl\{\birth(C'):C'\in \H,\, I(\basis(C'),\th)=1\Bigr\}
\end{equation}
with the convention $\min\emptyset = 0$. By definition, $\ti(\H,\Lambda,\th)
\le 0$.

\begin{thm}\label{thm2}
  The clan $\A^{\Lambda,0}$ is the limit as $t\to\infty$ of a process $\A_t$,
  defined by the initial condition $\A_0=\emptyset$ and the evolution equation
\begin{eqnarray}
  \label{oo1}
\lefteqn{
  \E\,\Bigl({d F(\A_t)\over dt}\Given \A_s,\, 0\le s\le t\Bigr) 
  \;=\;} \nonumber\\
&&\int_{\G(\A_t,\Lambda)} \nu(d\th) 
\int_{t+\ti(\A_t,\Lambda,\th)}^\infty   ds\,\,e^{-s}\,
   \Bigl[F(\A_t\cup (\th,-t,s)) - F(\A_t)\Bigr]\;.
\end{eqnarray}
Here $F$ is an arbitrary function depending on a finite number of individuals
intersecting $\Lambda$ and we have denoted $(\th,-t,s)$ the cylinder of base
$\th$, born at time $-t$ and with lifetime $s$.
\end{thm}

For completeness, a proof of this theorem is presented in the Appendix.  For
the free discrete loss-network (contours, random cluster) processes,
\begin{equation}
  \label{eq:ee1}
\int_{\G(\A_t,\Lambda)} \nu(d\th) \,F(\th)\; =\;
\sum_{\th\in\G(\A_t,\Lambda)} w(\th) F(\th)
\end{equation}
while for free birth-and-death processes on $\R^d$, 
\[
\int_{\G(\A_t,\Lambda)} \nu(d\th) \, F(\th)\;=\;
 \int_{\R^d} f(x)\, dx \int_{\G_x}\; \pi_x(dg) \;\one\{x\oplus g\in
\G(\A_t,\Lambda)\}\,F(x,g)\;.
\]

Notice that $\A_t$ is a monotone process ($\A_t\subset\A_{t+s}$) in which at
time $t$ only cylinders in $\G(\A_t,\Lambda)$ can be included. The inclusion
of a cylinder born at time $-t$ requires that either (a) its basis is
incompatible with that of some cylinder born later and its lifespan reaches
the birth-time of such cylinder, or (b) its basis is compatible with those of
all cylinders born later, but it intersects $\Lambda$ and the cylinder
survives up to time equal zero.  The last condition is ensured via the
convention $\min\emptyset=0$ in the definition of $\ti$.

\paragraph*{Algorithm to construct the backwards clan of a finite
  region} 

The combination of \reff{bpa}/\reff{oo1} can be translated into the following
explicit algorithm. We do it first for the case of countable number of
individuals and indicate at the end of this section how to proceed in the
continuous case. To generate $\A^{\Lambda,0}$:

\begin{enumerate}

\item Start with $\tau=0$ and $\H=\emptyset$.
  
\item Let $\H$ be the current set of cylinders and $\tau$ the current
  $-\min\{\birth(C): C\in \H\}$. For each $\ga\in \G(\H,\Lambda)$ generate an
  independent realization of the first time $\tau_1(\ga)$ of the
  non-homogeneous Poisson process in $\R$ with intensity
\begin{equation}
\label{la4}
\lambda_\ga(ds):=w(\ga)\,e^{-s+\ti(\H,\Lambda,\ga)}\,\one\{s>\tau\}\,ds \;. 
\end{equation}
Notice that $\tau_1(\ga)$ may be infinity.

\item Order the set $\{\tau_1(\ga):\ga\in\G(\H,\Lambda)\}$. Let
  $\widetilde\tau$ be the infimum of this set. This is well defined and
  strictly positive because the condition $\alpha<1$ guarantees that the total
  rate $\sum_{\ga\in \G(\H,\Lambda)} \int_{\R^+}\lambda_\ga(ds)<\infty$.
  
\item If $\widetilde\tau<\infty$, call $\ga_1$ the basis corresponding to the
  minimum (i.e.\/ $\tau_1(\ga_1) = \widetilde\tau$). Update
  $\tau\leftarrow\widetilde\tau$ and $\H \leftarrow
  \H\cup\{(\ga_1,-\tau,\tau+\ti(\H,\Lambda,\ga)+R_1)\}$, where $R_1$ is an
  exponential random variable with rate 1 independent of everything. In the
  sequel ignore the set $\{\tau_1(\ga)\,:\,I(\ga,\ga_1)=1\}$ (we can reuse the
  remaining $\tau_1$) and go to (2).
  
\item If $\widetilde\tau=\infty$ set $\A^{\Lambda,0}=\H$ and stop.  By Theorem
  \ref{thm2} the distribution of the set $\A^{\Lambda,0}$ so generated is
  exactly that of the free birth-and-death process.

\end{enumerate}

If $\tau_i$ are the successive times of jump of $\A_t$, then $\H_i =
\A_{\tau_i}$ have the same distribution as the $i$-th iterate of the above
algorithm.

In the continuous case, time and space cannot be in general separated.
Instead of steps (2) and (3) above we must consider the events $(\gamma,s)$ of
a Poisson process on $\G\times\R^+$ with intensity
\begin{equation}
  \label{cii}
  \lambda(d(\ga,s)) = \nu(d\ga)\,
  e^{-s+\ti(\H,\Lambda,\ga)}\,\one\{s>\tau\}\,\one\{\ga\in\G(\H,\Lambda)\}\,ds
  \;. 
\end{equation}
For a finite window $\Lambda$ the total rate is finite, hence these events can
be well ordered by looking to the time coordinate.  If the set of these events
is not empty, we take $\widetilde\tau$ to be the minimal time coordinate (it
is strictly positive with probability one) and denote $\ga_1$ the associated
individual.  If the Poisson process with rate density \reff{cii} yields no
event we take $\widetilde\tau=\infty$.  We then continue as in (4).

This algorithm plus the subsequent ``cleaning algorithm'' constitutes our
perfect simulation scheme.

\paragraph*{The cleaning algorithm} 

Let $\A^{\Lambda,0}$ be the clan of the cylinders whose life contains time 0
and basis intersects $\Lambda$. The following algorithm shows how to construct
inductively the set $\K^{\Lambda,0}$ of kept cylinders.

\begin{enumerate}
\item Start with $\H=\A^{\Lambda,0}$ and $\K=\emptyset$ ($\H$ is formed by the
  cylinders to be tested and $\K$ by those already kept).
  
\item If $\H$ is empty go to 5. If not, order the cylinders of $\H$ by time of
  birth. Let $C_1$ be the first of those cylinders; call $\ga_1$ its basis and
  $\tau_1$ its birth-time. Let $\xi_1$ be the set of bases of the cylinders in
  $\K$ alive at $\tau_1$ which are incompatible with the basis of $C_1$. Let
  $Z_1$ be a random variable uniformly distributed in $[0,1]$ independent of
  everything.
  
\item If $Z_1< M(\ga_1\vert\xi_1)$, then update: $\H \leftarrow
  \H\setminus\{C_1\}$, $\K \leftarrow \K\cup\{C_1\}$. Go to 2.
  
\item If $Z_1> M(\ga_1\vert\xi_1)$, then update: $\H \leftarrow
  \H\setminus\{C_1\}$. Go to 2.
  
\item Set $\K^{\Lambda,0} = \K$ and stop.  By Theorem \ref{thm1} (ii) the
  distribution of this clan $\K^{\Lambda,0}$ is exactly that of the
  interacting birth-and-death process.
\end{enumerate}

\paragraph*{Algorithm to simulate a finite window of $\mu$}

This is the easiest part. Once the set $\K^{\Lambda,0}$ of kept cylinders has
been determined, take the configuration $\eta$ defined by
\begin{equation}
  \label{set}
 \eta(\ga) = \sum_{C\in \K^{\Lambda,0}} \one\{C \hbox{ has basis }\ga
 \hbox{ and life containing }0\} 
\end{equation}
for $\ga$ intersecting $\Lambda$.  This configuration has the marginal
distribution of the infinite-volume measure $\mu$ on the (not necessarily
finite) set $\G_\Lambda = \N^{\{\th\in\G:\th\cap\Lambda\neq \emptyset\}}$.
This fact is guaranteed by Theorem \ref{thm1} (ii).

\section{Errors in Perfect simulation?}

Even in finite volume, perfect simulation algorithms are subjected to error.
In general terms, a perfect simulation algorithm of a measure $\mu$ on a set
$\X$ is a function $\Phi:[0,1]^\N\to \X$, such that, if $(U_n)_{n\in\N}$ is a
sequence of i.i.d. uniform in $[0,1]$ random variables, there exists a
stopping time $T$ for $(U_n)$ such that $\Phi$ depends only on the first $T$
coordinates of $(U_1,U_2,\dots)$ and
\begin{equation}
  \label{psa}
  \P(\Phi(U_1,\dots,U_T) \in A) = \mu(A)\;.
\end{equation}
The CFTP algorithm, for instance, stops when a random value $t$ is found such
that the different copies of the algorithm coupled from time $-t$ started with
all possible initial conditions attain the same configuration at time $0$.
Finding $t$ requires the use of a random number $T(t)$ of uniform random
variables, which must be less than $S=$``the maximum time left in order to
have the results ready for the next congress'', for instance.  Thus, one
actually samples from the distribution defined by
\begin{equation}
  \label{ps2}
  \P(\Phi(U_1,\dots,U_T) \in A \,|\, T<S)
\end{equation}
which is different from, though as $S\to\infty$ converges to, \reff{psa}.
This is the so-called {\em impatient-user bias}. The CFTP algorithm also
permits the construction of a joint realization $(\eta,\xi)$ with marginals
\reff{psa} and \reff{ps2} such that $T<S$ implies $\eta=\xi$. In fact, as
pointed in Proposition 6.2 of Fill (1998)
\begin{equation}
  \label{b:ps3}
  \sup_{A} \Bigl\vert \P\Bigl(\Phi(U_1,\dots,U_T) \in A \,\Bigm|\, 
T<S\Bigr) - \mu(A) \Bigr  \vert \,\leq\,
 \frac{\P[T>S]}{1-\P[T>S]}\;.
\end{equation}
In our algorithm $T$ is determined by the number of uniform random variables
necessary to construct the clan of the observed region $\Lambda$.

When the possible sizes of the individuals $\gamma$ form an unbounded set, for
instance for the Peierls contours of the Ising model, practical limitations
prevent the inclusion of all possible sizes in the simulation.  In fact the
mere enumeration of the possible contours is beyond reach when more than a few
dozens of links are involved.  This is tantamount to a ``space impatient-user
bias'': the user is forced to do a space cut-off that produces a bias, even
when the actual probability for a cut event to take place is tiny.  In
mathematical terms, one actually samples from the conditioned measure
\begin{equation}
  \label{ps4}
  \P\Bigl(\Phi(U_1,\dots,U_T) \in A \,\Bigm|\, \{K<k\}\cap\{T<S\}\Bigr)
\end{equation}
where $K=$ ``maximum perimeter of bases of cylinders in the clan'' ($k=30$,
for instance).  In fact, our approach also admits a joint realization
$(\eta,\xi)$ with marginal distributions \reff{psa} and \reff{ps4} such that
$\eta=\xi$ if $K<k$ and $T<S$, and such that $\P(\{K\ge k\}\cap\{T>S\})$ goes
to zero exponentially fast in $S$ and in the cutoff of the length of the
contours (30 in our example).  Slightly more precisely, a bound like
\reff{b:ps3} holds with
\begin{equation}
  \label{eq:pss5}
  \P\Bigl(\{K\ge k\}\cap\{T>S\}\Bigr) \; \le \; 
{\rm O}\Bigl(\alpha^T\, \times \, \sup_x\pi_x(K>k)\Bigr)\;.
\end{equation}
This follows from the subcriticality of the majorizing branching process.  For
the Ising model, for instance, $\pi_x(K>k)={\rm O}(e^{-\beta k})$.

\section{Conclusion}

Our algorithm offers an approach to perfect simulations of processes with
infinite state space.  The fact that there is no coupling between different
initial conditions, makes it a flexible tool for processes with a large state
space. No ``sandwiching processes'' need to be followed; the free process is a
natural ``dominating process'' in our setting.  In addition, our algorithm is
backed by a rather detailed theory that allows the estimation of various
properties of the resulting measure, as well as possible errors.  In
particular, our approach is not free from the ``impatient-user bias'', but the
resulting error is relatively straightforward to control.

A noteworthy feature of our approach is that the perfect simulation stage
applies, in fact, to the free process.  Interacting processes are then
obtained by a deterministic ``cleaning''.  As a consequence our scheme allows
the simultaneous simulation of all processes absolutely continuous with
respect to the same free process.  This coupled construction could be
potentially useful, for instance to establish comparison criteria.

The algorithm admits a further generalization more or less immediate that has
not been pursued here: it can be applied to processes with variable death rate
that, however, must be uniformly bounded from below by 1, say. The dependences
in the birth and death-rates induce definitions of incompatibility and
respective parameters $\alpha$ (cf.  \reff{t54a}, \reff{mn2}). A construction
analogous to the one described in this paper can be performed but with a
thinning algorithm that takes also into account the variable death-rates.

In this work, the advantages of the approach have been exploited only at a
theoretical level, where it has led to a new treatment of systems with
exclusions and to better estimates of regions of existence of a number of
processes. Berthelsen and M{\o}ller (2001) compared it to the dominated CFTP
introduced by Kendall and M{\o}ller (2000). Based on simulation results, the
authors show that the dominated CFTP is better than the algorithm based on the
clan of ancestors in the particular case of a Strauss process (see Equation
\reff{eq:300ca}) defined on a unit square with $e^\beta_1 = 100$ and
$e^\beta_2 = 0$ (the so-called hard-core process), $0.5$ and $1$ (a Poisson
processes with rate 100). This is obviously the case from the description of
the processes since the backward construction of our algorithm stops when the
dominated Poisson process regenerates and usually the coupling of CFTP is
achieved before it in the finite case. However, it should be noticed that the
algorithm based on the clan of ancestors was designed for sampling the
infinite-volume process viewed in a finite window.  This seems to be a much
more interesting and challenging problem which has been studied by M\'{a}ric
(2002) for the specific case of one-dimensional loss networks with bounded
calls. No comparison was made to other perfect simulation schemes.

Finally, we hope that a suitable combination of our ideas with some rejection
sampling scheme could yield a version free of the user-impatience bias.

\appendix
\section{Proof of Theorems}
\def\bS{{\bf S}}
\paragraph*{ Proof of Theorem \protect\ref{thm1}}

We need to show that $\eta^{\zeta}_t$ has generator \reff{genp}.  Denote
$\eta_t = \eta^\zeta_t$ and $\K[0,t]$ the set of kept cylinders born at time
zero or after time zero (this includes the cylinders induced by the initial
configuration) and for $F$ a function depending on individuals intersecting a
region with finite total rate, write
\begin{eqnarray}
  \label{153}
  \lefteqn{ [F(\eta_{t+h}) -  F(\eta_t)] }\nonumber\\
  &=& \sum_{C\in \K[0,t+h]} \one\{\birth(C)\in [t,t+h]\}
  [F(\eta_t+\de_{\basis(C)})-F(\eta_t)]\nonumber\\
  &&{}+\sum_{C\in \K[0,t]} \one\{\life(C)\ni t,\, \life(C)\not\ni
  t+h\}[F(\eta_t-\de_{\basis(C)})-F(\eta_t)] \nonumber\\
  &&{}+ \{\hbox{other things}\},
\end{eqnarray}
where $\{${other things}$\}$ refer to events with more than one Poisson mark
in the time interval $[t,t+h]$ for the contours in the (finite) support of
$F$. Since the total rate of the Poisson marks in this set is finite, the
event $\{$other things$\}$ has a probability of order $(h\,m(\supp(F),\G))^2$,
where $m$ is defined in \reff{mmm}.
We have
\begin{eqnarray}
  \label{154} \lefteqn{ \sum_{C\in \C} \one\{\birth(C)\in
  [t,t+h]\}\,\one\{C\in
  \K[0,t+h]\}\,[F(\eta_t+\de_{\basis(C)})-F(\eta_t)]}\nonumber\\ &
  =&\sum_{C\in \C} \one\{\birth(C)\in
  [t,t+h]\}\,\one \{\flag(C)<M(\basis(C)\vert \eta_t)\}\nonumber\\ 
&&\qquad\qquad\qquad\qquad\qquad
  \times[F(\eta_t+\de_{\basis(C)})-F(\eta_t)]
\end{eqnarray}

To compute the second term of \reff{153}, observe that $\life(C)$ is
independent of $\birth(C)$ and both the event $\{C\in \K[0,t]\}$ and $\eta_t$
are ${\mathcal F}_t$-measurable. Here ${\mathcal F}_t$ is the $\sigma$-algebra
generated by the births and deaths occurred before $t$.  Hence
\begin{eqnarray}
\lefteqn{\P\Bigl(\life(C)\ni
  t,\, \life(C)\not\ni
  t+h \Given{\mathcal F}_t\Bigr)}\nonumber\\
  &=& \P\Bigl(\life(C)\not\ni t+h \Given \life(C)\ni
  t\Bigr)\,\one\{\life(C)\ni t\} 
 \nonumber
\end{eqnarray}
and
\begin{eqnarray}
  \label{155}
 \lefteqn{
 \E \Bigl[ \sum_{C} \one\{C\in\K[0,t]\}\,\one\{\life(C)\ni
  t,\, \life(C)\not\ni
  t+h\}}\nonumber\\ 
&&\qquad\qquad\qquad\qquad\qquad\times\,
[F(\eta_t-\de_{\basis(C)})-F(\eta_t)]\Bigr]  
\nonumber \\
 & =&
  \E \Bigl[  \sum_C \P\Bigl(\life(C)\not\ni
  t+h\Given \life(C)\ni t\Bigr)\,\one\{C\in\K[0,t],\;\life(C)\ni
  t\}\,\nonumber \\
&&\qquad\qquad\qquad\qquad\qquad\times\,
  [F(\eta_t -\de_{\basis(C)})-F(\eta_t)]\Bigr].
\end{eqnarray}
Since $\life(C)$ is exponentially distributed with mean 1,
\begin{eqnarray}
\P\Bigl(\life(C)\not\ni t+h \Given \life(C)\ni t\Bigr) 
&=&  h + o(h)\,. \label{1554}
\end{eqnarray}
Taking the expectation of \reff{153} and substituting
\reff{154}--\reff{1554} we get
\begin{eqnarray}
  \label{156} \lefteqn{\E [F(\eta_{t+h}) - F(\eta_t)] }\nonumber\\ &
  = &h\, \int_\G \nu(d\ga) \,\E\Bigl(M(\ga\vert \eta_t)\,
    [F(\eta_t+\de_{\ga})-F(\eta_t)]\Bigr) + o(h)\nonumber \\ 
& & \qquad\qquad+ \;h\,\sum_{\ga:\eta_t(ga)>0}\,\E \Bigl(\eta_t(\ga)
\,[F(\eta_t-\de_{\ga})-F(\eta_t)] \Bigr)+o(h) 
\end{eqnarray}
which dividing by $h$ and taking limit gives 
\begin{equation}
  \label{dge}
  {d\E F(\eta^\zeta_t)\over dt} \;=\; A \E F(\eta^\zeta_t). \qquad\qquad \square
\end{equation}

\paragraph*{Proof of Theorem \protect\ref{thm2}}

Define
\begin{equation}
  \label{65}
  \A_t\; =\; \Bigl\{C'\in \A^{\Lambda,0}:
  0>\birth(C')>-t\Bigr\}\; =\; \A^{\Lambda,0}\cap \C[-t,0] ,
\end{equation}
that is, the set of cylinders in $\A^{\Lambda,0}$ with birth-time posterior to
$-t$. It suffices to prove that the process so defined satisfies the evolution
equation \reff{oo1}.

The inclusion of a new cylinder in the time interval $[t,t+h]$ depends on the
existence of a birth Poisson mark in $[-t-h,-t]$ whose corresponding cylinder
is incompatible with some $C'\in \A_t$.  That is, if $C$ is a cylinder with
$I(\basis(C'),\basis(C))=1$ for some $C'\in \A_t$,
\begin{eqnarray}
&&\P \Bigl(\A_{t+h}=\widetilde\A\cup
 C\Given\A_t=\widetilde \A,\, \A_{t'}=\widetilde
 \A_{t'}, \,t'\in[0,t) \Bigr)\;=\; \P \Bigl\{C\in\C\,:\, \nonumber\\ 
 && \birth(C)\in [-t-h,-t]\,,\,
 \death(C)>t-\ti(\widetilde\A,\Lambda,\basis(C))\Bigr\}\;+ \;o(h)  \nonumber
\end{eqnarray}
The remainder $o(h)$ is the correction related to the probability that $C$ is
not the only \emph{relevant} cylinder born in $[-t-h,-t]$.  Hence
\[
o(h)\le \Bigl(h\,\sum_{C\in \A_t} m(\basis(C),\G)\Bigr)^2\,\le
h^2\, |\A_t|^2 \,\alpha^2
\]
where $m$ is defined in \reff{mmm} and $|\A_t|$ stands for $\sum_{C\in\A_t}
q(\basis(C))$, where $q$ is the measure used to define $\alpha$ in
\reff{t54a}.  Since the birth-time is independent of the lifetime which is
exponentially distributed with rate one,
\begin{eqnarray}
 \lefteqn{\P \Bigl(\A_{t+h}=\widetilde\A\cup
 C\Given\A_t=\widetilde \A ,\, \A_{t'}=\widetilde
 \A_{t'}, \,t'\in[0,t)  \Bigr) }\nonumber\\ 
 &=& \P \Bigl\{C\in\C\,: \,\birth(C)\in [-t-h,-t]\Bigr\}\,\nonumber\\ 
&&\qquad\times
\P\Bigl( \life(C)>t-\ti(\widetilde\A,\Lambda,\basis(C))\Bigr)\;
 +\; o(h)  \nonumber\\
 &=& h\,f(\basis(C))\, e^{-t+\ti(\widetilde\A,\Lambda,\basis(C))} 
\;+\; o(h)\;.\label{oo2}
\end{eqnarray}
This implies that when the configuration at time $t^-$ is $\widetilde\A$, a
new cylinder with basis $\gamma$ is included in $\A_t(\Upsilon)$ at rate
\begin{equation}
  \label{ooo}
  f(\gamma)\,\, e^{-t+\ti(\widetilde\A,\Lambda,\gamma)}.
\end{equation}
From \reff{oo2}, as in the computation of the forward Kolmogorov equations, we
get \reff{oo1}.  This equation characterizes the law of the process
$\A_t(\Upsilon)$ as a non-homogeneous Markov process.\hfill\square

\section*{Acknowledgments} It is a pleasure to thank W.~Kendall,
J.~M\o ller and E.~Th\"onnes for enlightening comments.  We also thank Timo
Seppalainen for a nice discussion on the random cluster model. We thank two
referees for their comments and criticisms that helped to improve the paper.
This paper was written while PAF was visiting professor at the Laboratoire des
Probabilit\'es de l'Universit\'e de Paris VI, the Departement de Mathematiques
de l'Universit\'e de Cergy Pontoise and the UPRES-A CNRS 6085 de
l'Universit\'e de Rouen.

This work was partially supported by FAPESP, CNPq and FINEP (N\'ucleo de
Excel\^encia ``Fen\^omenos cr\'\i ticos em probabilidade e processos
estoc\'asticos'' PRONEX-177/96).

\newpage

\parskip 0pt

\obeylines
Roberto Fern\'andez, 
UPRES-A, CNRS 6085, 
Math\'ematiques, Site Colbert, 
Universit\'e de Rouen, 
F 76821 Mont Saint Aignan -- Cedex, 
FRANCE 
{\tt Roberto.Fernandez@univ-rouen.fr} 
http://www.univ-rouen.fr/upresa6085/Persopage/Fernandez/ 
\vskip 5mm

Pablo A. Ferrari
IME USP, 
Caixa Postal 66281, 
05311-970 - S\~{a}o Paulo,
BRAZIL
{\tt pablo@ime.usp.br}
http://www.ime.usp.br/\~{}pablo
\vskip 5mm

Nancy L. Garcia
IMECC, UNICAMP, Caixa Postal 6065, 
13081-970 - Campinas SP 
BRAZIL
{\tt nancy@ime.unicamp.br}
http://www.ime.unicamp.br/\~{}nancy


\begin{thebibliography}{FMMR99}


\bibitem{badlie95}
A.~J. Baddeley and M.~N.~M. van Lieshout (1995).
\newblock Area-interaction point processes.
\newblock {\em Ann. Inst. Statist. Math.}, 47(4):601--619.

\bibitem{bm} Kasper K. Berthelsen and Jesper Møller (2001). Spatial jump
  processes and perfect simulation. Preprint.


\bibitem{vdb93} J. van den Berg (1993).
\newblock  A uniqueness condition for
  Gibbs measures, with application to the $2$-dimensional Ising
  antiferromagnet. Comm. Math. Phys. 152 (1993), no. 1, 161--166. 
  
\bibitem{vdBMae94} J. van den Berg and C. Maes (1994).
\newblock Disagreement
  percolation in the study of Markov fields. {\em Ann.  Probab. \bf
    22}(2):749--763.
  
\bibitem{vdBSte99} J. van den Berg and J. E. Steif (1999).
  \newblock On the existence and nonexistence of finitary codings for
  a class of random fields.  {\sl Ann. Probab. \bf 27}(3):1501--1522.
  
\bibitem{bry86} D.~C. Brydges (1986).  \newblock A short course
  on cluster expansions.  \newblock In {\em Ph\'enom\`enes critiques,
    syst\`emes al\'eatoires, th\'eories de gauge, Part I, II (Les
    Houches, 1984)}, pages 129--183. North-Holland, Amsterdam-New
  York.

\bibitem{ct01} H.~Cai (1999)  \newblock Exact sampling using
  auxiliary variables.  \newblock Preprint.

\bibitem{cff02} F.~ Comets, R.~ Fernandez and
  P.~A.~Ferrari (2002). Processes with Long Memory: Regenerative
  Construction and Perfect Simulation. To appear \newblock {\em
    Ann. Appl. Probab.}, arXiv:math.PR/0009204

\bibitem{ct01} J.~N.~Corcoran and R.~ L.~Tweedie.  \newblock
  Perfect sampling of ergodic Harris chains (2001). \newblock {\em
    Ann. Appl. Probab.}, {\bf 11}(2): 438--451.

\bibitem{dob96} R.~L. Dobrushin (1996).  \newblock Perturbation
  methods of the theory of {G}ibbsian fields.  \newblock In {\em
    Lectures on probability theory and statistics (Saint-Flour,
    1994)}, pages 1--66. Springer, Berlin.
  
\bibitem{fer89} P. A. Ferrari (1990). Ergodicity for spin systems with
  stirrings. {\sl Ann. Probab. \bf 18}(4):1523--1538.
  
\bibitem{ffg98} R.~Fern{\'a}ndez, P.~A. Ferrari, and N.~L.
  Garcia (1998).  \newblock Measures on contour, polymer or animal
  models. {A} probabilistic approach.  \newblock {\em Markov Process.
    Related Fields}, 4(4):479--497.  \newblock I Brazilian
  School in Probability (Rio de Janeiro, 1997).
  
\bibitem{ffg01} R.~Fern{\'a}ndez, P.~A. Ferrari, and N.~L.
  Garcia (2001).  \newblock Loss network representation of {P}eierls
  contours.  \newblock {\em Ann. Probab.}, {\bf 29}(2): 902--937
  
\bibitem{ferfrosok92} R. Fern{\'a}ndez, J. Fr{\"o}hlich and
  A. D. Sokal (1992). \newblock {\sl Random Walks, Critical Phenomena,
  and Triviality in Quantum Field Theory}, Springer-Verlag,
  Berlin--Heidelberg--New York.


\bibitem{fg98} P.~A. Ferrari and N.~L. Garcia (1998).  \newblock
  One-dimensional loss networks and conditioned ${M}/{G}/\infty$ queues.
  \newblock {\em J. Appl. Probab.}, 35(4):963--975.
  
\bibitem{fil98} J.~A. Fill(1998).  \newblock An interruptible
  algorithm for perfect sampling via {M}arkov chains.  \newblock {\em
    Ann. Appl. Probab.}, 8(1):131--162.

\bibitem{filhub00} J.~A. Fill and M.~Huber (2000).
  \newblock The Randomness Recycler: A New Technique for Perfect
  Sampling.  \newblock Preprint. arXiv:math.PR/0009242. 
   
  
\bibitem{filetal99} J.~A. Fill, M.~Machida, D.~J. Murdoch, and
  J.~S.  Rosenthal (2000).  \newblock Extension of {F}ill's perfect
  rejection sampling algorithm to general chains. Proceedings of the
  Ninth International Conference "Random Structures and Algorithms"
  (Poznan, 1999). {\sl Random Structures and Algorithms.} no. 3-4, 290--316.   

\bibitem{fk72} C.~M. Fortuin, P.W. Kasteleyn (1972) On the random cluster
model. I. Introduction and relation with other models. {\sl Physica \bf 57}
536--564. 

\bibitem{fostwe98} S.~G. Foss and R.~L. Tweedie (1998). \newblock
  Perfect simulation and backward coupling. {\sl Stoch. Models}, {\bf
  14}, pp. 187--203.

\bibitem{gar00} N.L. Garcia (2000). Perfect simulation of spatial 
processes, {\it Resenhas IME-USP}, {\bf 4}(3), 281--324.

\bibitem{gawkotkup87} K. Gaw\c edzki, R. Koteck{\`y} and A.
  Kupiainen (1987).  Coarse graining approach to first-order phase
  transitions.  {\sl J. Stat. Phys.}, {\bf 47}, pp\ 701--724.
  
\bibitem{grimmett} G.~Grimmett (1995).  \newblock The
  stochastic random-cluster process and the uniqueness of
  random-cluster measures.  \newblock {\em Ann. Probab.},
  23(4):1461--1510.
  
\bibitem{hagliemol96} O.~{H{\"a}ggstr{\"o}m}, M.~N.~M. van
  Lieshout, and J.~M{\o}ller (1999).  \newblock Characterization
  results and {M}arkov chain {M}onte {C}arlo algorithms including
  exact simulation for some spatial point processes.
\newblock  {\em Bernoulli}, {\bf 5}(4):641--658.

\bibitem{HasSte99} O.~{H{\"a}ggstr{\"o}m}, J. Steif (2000).
  \newblock Propp--Wilson algorithms and finitary codings for high
  noise random fields.  {\sl Combin. Probab. Comput.}{\bf 9}(5): 425--439.


\bibitem{hc00} (2000). J.~ P.~ Hobert and C.~ P.~
    Robert (2000). \newblock Moralizing perfect sampling. Preprint.

\bibitem{Kelly}Kelly, F. P. (1991) Loss networks.\newblock 
\newblock {\em Ann. Appl. Probab. \bf 1} 3:319--378. 


\bibitem{KelRip76} F. P. Kelly and B. D Ripley (1976).
\newblock A note on the Strauss' model for clustering
\newblock {\em Biometrika \bf 63}, 357--360.


\bibitem{ken97} W.~S. Kendall (1997).  \newblock On some
  weighted {B}oolean models.  \newblock In D.~Jeulin, editor, {\em
    Proceedings of the International Symposium on Advances in Theory
    and Applications of Random Sets (Fontainebleau, 1996)}, pages
  105--120. World Sci. Publishing, River Edge, NJ.
  
\bibitem{ken98} W.~S. Kendall (1998).  \newblock Perfect
  simulation for the area-interaction point process.  \newblock In
  L.~Accardi and C.~C. Heyde, editors, {\em Probability Towards 2000},
  pages 218--234. Springer.

\bibitem{kenmol99} W.~S. Kendall and J.~M{\o}ller (2000). \newblock
  Perfect simulation using dominating processes on ordered spaces, with
  applications to locally stable point processes, \newblock {\sl Adv. Appl.
    Probab.} {\bf 32}(3): 844--865.   
  
\bibitem{kotpre86} R.~Koteck{\'y} and D.~Preiss (1986).
  \newblock Cluster expansion for abstract polymer models.  \newblock
  {\em Comm. Math. Phys.}, 103(3):491--498.
  
\bibitem{ma02} N. M\'{a}ric (2002) \newblock Perfect simulation for a
  continuous one-dimensional loss network. \newblock Master's thesis,
  IMECC/UNICAMP. Available at
  http://www.ime.unicamp.br/rel$\underline{\phantom o}$pesq/2002/rp18-02.html

  
\bibitem{Moe00} J. M\o ller (2000).\newblock A review on
  perfect simulation in stochastic geometry {\sl IMS Lecture Notes,
    Monograph Series}.
  
\bibitem{olipic90} E.\ Olivieri and P.\ Picco (1990).  \newblock
  Cluster expansion for $D$-dimensional lattice systems and finite
  volume factorization properties. \newblock {\sl J. Stat. Phys.},
  {\bf 59}, pp.\ 221--256.

\bibitem{prowil96}
J.~G. Propp and D.~B. Wilson (1996).
\newblock Exact sampling with coupled {M}arkov chains and applications to
  statistical mechanics.
\newblock In {\em Proceedings of the Seventh International Conference on Random
  Structures and Algorithms (Atlanta, GA, 1995)}, volume~9, pages 223--252.
  
\bibitem{str75} David~J. Strauss (1975).  \newblock A model for
  clustering.  \newblock {\em Biometrika}, 62(2):467--475.
  
\bibitem{tho97} E.~{Th{\"o}nnes} (1999).  \newblock
  Perfect simulation of some point processes for the impatient user.
  \newblock {\em Adv. Appl. Probab. \bf 31} 69--87.

\end{thebibliography}
\end{document}